\documentclass[10pt]{article}%
\usepackage{amsmath}
\usepackage{amsfonts}
\usepackage{amssymb}
\usepackage{graphicx}
\usepackage{color}%
\providecommand{\U}[1]{\protect\rule{.1in}{.1in}}
\begin{document}

\title{Observation estimate for the heat equations \\with Neumann boundary condition \\via logarithmic convexity}
\author{R\'{e}mi Buffe\thanks{ Institut Elie Cartan de Lorraine Universit\'{e} de
Lorraine, Site de Nancy \& Inria (Project-Team SPHINX) B.P. 70239, F-54506 Vandoeuvre-l\`{e}s-Nancy Cedex,
France. E-mail address: remi.buffe@univ-lorraine.fr} , Kim Dang Phung\thanks{
Institut Denis Poisson, Universit\'{e} d'Orl\'{e}ans, Universit\'{e} de Tours
\& CNRS UMR 7013, B\^{a}timent de Math\'{e}matiques, Rue de Chartres, BP.
6759, 45067 Orl\'{e}ans, France. E-mail address: kim\_dang\_phung@yahoo.fr}}
\date{}
\maketitle

\bigskip

Abstract .- We prove an inequality of H\"{o}lder type traducing the unique
continuation property at one time for the heat equation with a potential and
Neumann boundary condition. The main feature of the proof is to overcome the
propagation of smallness by a global approach using a refined parabolic
frequency function method. It relies with a Carleman commutator estimate to
obtain the logarithmic convexity property of the frequency function.

\bigskip


\bigskip

\section{Introduction and main result}

\bigskip

In this paper, we establish the observation inequality at one time for the
heat equation with a potential and Neumann boundary condition. The analysis is
based on the parabolic frequency function method \cite{K} adjusted for a
global approach.

\bigskip

Let $\Omega\subset\mathbb{R}^{\text{n}}$ be a bounded connected open set with
boundary $\partial\Omega$ of class $C^{\infty}$. Consider in $\{(x,t)\in
\Omega\times\left(  0,T\right)  \}$ the heat equation with a potential and
Neumann boundary condition
\[
\left\{
\begin{array}
[c]{ll}%
{\partial}_{t}u-\Delta u+au=0\text{ ,} & \quad\text{in}~\Omega\times\left(
0,T\right)  \text{ ,}\\
\partial_{n}u=0\text{ ,} & \quad\text{on}~\partial\Omega\times\left(
0,T\right)  \text{ ,}\\
u\left(  \cdot,0\right)  \in L^{2}(\Omega)\text{ .} &
\end{array}
\right.
\]
Here, $T>0$, $a\in L^{\infty}\left(  \Omega\times\left(  0,T\right)  \right)
$ and $n$ is the unit outward normal vector to $\partial\Omega$.

\bigskip

We propose the following result.

\bigskip

Theorem 1 .- \textit{Let }$\omega$\textit{ be a non-empty open subset of
}$\Omega$\textit{. For any }$t\in\left(  0,T\right]  $\textit{, }%
\[
\left\Vert u\left(  \cdot,t\right)  \right\Vert _{L^{2}\left(  \Omega\right)
}\leq\left(  e^{K\left(  1+\frac{1}{t}+t\left\Vert a\right\Vert _{L^{\infty
}\left(  \Omega\times\left(  0,t\right)  \right)  }+\left\Vert a\right\Vert
_{L^{\infty}\left(  \Omega\times\left(  0,t\right)  \right)  }^{2/3}\right)
}\left\Vert u\left(  \cdot,t\right)  \right\Vert _{L^{2}\left(  \omega\right)
}\right)  ^{\beta}\left\Vert u\left(  \cdot,0\right)  \right\Vert
_{L^{2}\left(  \Omega\right)  }^{1-\beta}\text{ .}%
\]
\textit{Here }$K>0$\textit{ and }$\beta\in\left(  0,1\right)  $\textit{ only
depend on }$\left(  \Omega,\omega\right)  $\textit{.}

\bigskip

Such observation estimate traduces the unique continuation property at one
point in time saying that if $u=0$ in $\omega\times\left\{  t\right\}  $, then
$u$ is identically null. Applications to bang-bang control and finite time
stabilization are described in \cite{PWZ} and \cite{BuP}. Our result is an
interpolation estimate which is more often used in a local way with a
propagation of smallness procedure (\cite{AEWZ}, \cite{FV}). Here the way we
choose to establish our main theorem is based on a global approach.

\bigskip

Recall that Theorem 1 implies the observability estimate for the heat
equations with a potential and Neumann boundary condition \cite{PW}. It is
well-known that the observability estimate for the heat equations can be
obtained from Carleman inequalities. In the literature, at least two
approaches allow to derive Carleman inequalities for parabolic equations: A
local one based on the Garding inequality and interpolation estimates for the
elliptic equations (\cite{LR}, \cite{LRL}, \cite{BM}); A global one based on
Morse functions and integrations by parts over $\Omega\times\left(
0,T\right)  $ (\cite{FI}, \cite{FGGP}). Besides, unique continuation results
can be deduced either by Carleman techniques or by logarithmic convexity of a
frequency function \cite{EFV}. Here we construct a new frequency function
adapted to the global approach. Further, we explicitly give the dependence of
the constants with respect to $\left\Vert a\right\Vert _{L^{\infty}}$ as in
\cite{FGGP}, \cite{DZZ}.

\bigskip

\bigskip

\section{Preliminaries}

\bigskip

In this section we derive three propositions on which our later results will
be based.

\bigskip

Proposition 1 .- \textit{Let }$\Omega\subset\mathbb{R}^{\text{n}}$\textit{ be
a bounded connected open set of class }$C^{\infty}$\textit{, and let }$\omega
$\textit{ be a non-empty open subset of }$\Omega$\textit{. Then there exist
}$d\in\mathbb{N}^{\ast}$\textit{, }$\left(  p_{1},p_{2},\cdot\cdot
,p_{d}\right)  \in\omega^{d}$\textit{ and }$\left(  \psi_{1},\psi_{2}%
,\cdot\cdot,\psi_{d}\right)  \in\left(  C^{\infty}\left(  \overline{\Omega
}\right)  \right)  ^{d}$\textit{ such that for all }$i\in\left\{  1,\cdot
\cdot,d\right\}  $

\begin{description}
\item[$\left(  i\right)  $] $\psi_{i}>0$\textit{ in }$\Omega$\textit{, }%
$\psi_{i}=0$\textit{ on }$\partial\Omega$\textit{ ,}

\item[$\left(  ii\right)  $] \textit{the critical points of }$\psi_{i}%
$\textit{ are nondegenerate,}

\item[$\left(  iii\right)  $] $\left\{  x\in\Omega;\left\vert \nabla\psi
_{i}\left(  x\right)  \right\vert =0\right\}  =\left\{  p_{j};j=1,\cdot
\cdot,d\right\}  $\textit{,}

\item[$\left(  iv\right)  $] $p_{i}$\textit{ is the unique global maximum of
}$\psi_{i}$\textit{ ,}

\item[$\left(  v\right)  $] \textit{for any }$j\in\left\{  1,\cdot
\cdot,d\right\}  $\textit{, }$\underset{\overline{\Omega}}{\text{max}}\psi
_{j}=\underset{\overline{\Omega}}{\text{max}}\psi_{i}$\textit{ .}
\end{description}

\bigskip

Remark .- $\left(  i\right)  $ implies that $\partial_{n}\psi_{i}\leq0$;
$\left(  iii\right)  $ says that the criticals points of $\psi_{i}$ are
isolated and form a discrete set; $\left(  iii\right)  $ implies that
$d=\sharp\left\{  x\in\Omega;\left\vert \nabla\psi_{i}\left(  x\right)
\right\vert =0\right\}  $\ and $\left\{  x\in\Omega;\left\vert \nabla\psi
_{i}\left(  x\right)  \right\vert =0\right\}  \subset\omega$.

\bigskip

Proof .- The existence of Morse functions (that is $C^{\infty}$ functions
whose critical points are nondegenerate) which are positive in $\Omega$ and
null on the boundary $\partial\Omega$ can be proved by virtue of the theorem
on the density of Morse functions (\cite[page 20]{FI}, \cite[page 80]{C},
\cite[Chapter 14]{TW}, \cite[page 433]{WW}). Next, by a small perturbation in
a small neighborhood of each critical points, no two critical points share the
same function value \cite[Theorem 2.34]{M}. Denote by $\psi$ such a smooth
function and let $a_{1},\cdot\cdot,a_{d}$ be its critical points such that
$\left\{  x\in\Omega;\left\vert \nabla\psi\left(  x\right)  \right\vert
=0\right\}  =\left\{  a_{j};j=1,\cdot\cdot,d\right\}  \subset\Omega$ and
$\psi\left(  a_{1}\right)  >\psi\left(  a_{2}\right)  >\cdot\cdot>\psi\left(
a_{d}\right)  $. Now we will move the critical points following the procedure
in \cite[Lemma 2.68]{C}. Introduce $p_{1},\cdot\cdot,p_{d}$ $d$\ points in
$\omega$ such that for each $i=1,\cdot\cdot,d$, there exists $\gamma_{i,j}\in
C^{\infty}([0,1];\Omega)$ be such that

\begin{itemize}
\item $\gamma_{i,j}$ is one to one for every $j\in\left\{  1,\cdot
\cdot,d\right\}  $,

\item $\gamma_{i,j}\left(  \left[  0,1\right]  \right)  \cap\gamma
_{i,l}\left(  \left[  0,1\right]  \right)  =\emptyset,\forall\left(
j,l\right)  \in\left\{  1,\cdot\cdot,d\right\}  $ such that $j\neq l$,

\item $\gamma_{i,j}\left(  0\right)  =a_{j}$, $\forall j\in\left\{
1,\cdot\cdot,d\right\}  $,

\item $\gamma_{i,j}\left(  1\right)  =\tau^{i-1}\left(  p_{j}\right)  $,
$\forall j\in\left\{  1,\cdot\cdot,d\right\}  $.
\end{itemize}

Here $\tau$ is $d$-cycle, that is $\tau\left(  p_{j}\right)  =p_{j+1}$ if
$j<d$\ and $\tau\left(  p_{d}\right)  =p_{1}$, $\tau^{0}=id$, $\tau^{i}%
=\tau^{i-1}\circ\tau$.

Introduce a vector field $V_{i}\in C^{\infty}(\mathbb{R}^{\text{n}}%
;\mathbb{R}^{\text{n}})$ such that $\overline{\left\{  x\in\mathbb{R}%
^{\text{n}};V_{i}\left(  x\right)  \neq0\right\}  }\subset\Omega$ and
$V_{i}(\gamma_{i,j}(t))=\gamma_{i,j}^{\prime}(t)$, $\forall j\in\left\{
1,\cdot\cdot,d\right\}  $. Let $\Lambda_{i}$ denote the flow associated to
$V_{i}$, that is $\partial_{t}\Lambda_{i}\left(  t,x\right)  =V_{i}\left(
\Lambda_{i}\left(  t,x\right)  \right)  $ and $\Lambda_{i}\left(  0,x\right)
=x$. One has $\Lambda_{i}\left(  0,a_{j}\right)  =a_{j}$, $\Lambda_{i}\left(
t,a_{j}\right)  =\gamma_{i,j}\left(  t\right)  $ and $\Lambda_{i}\left(
1,a_{j}\right)  =\tau^{i-1}\left(  p_{j}\right)  $. Further, for every
$t\in\mathbb{R}$, $\Lambda_{i}\left(  t,\cdot\right)  $ is a diffeomorphism on
$\Omega$ and $\Lambda_{i}\left(  t,\cdot\right)  \left\vert _{\partial\Omega
}\right.  =Id$. In particular, $\left(  \Lambda_{i}\left(  1,\cdot\right)
\right)  ^{-1}\left(  \tau^{i-1}\left(  p_{j}\right)  \right)  =a_{j}$.

It remains to check that $\psi_{i}:\overline{\Omega}\rightarrow\mathbb{R}$
given by $\psi_{i}\left(  x\right)  =\psi\left(  \left(  \Lambda_{i}\left(
1,\cdot\right)  \right)  ^{-1}\left(  x\right)  \right)  $ satisfies all the
required properties. Clearly, $\psi_{i}>0$ in $\Omega$, $\psi_{i}=0$ on
$\partial\Omega$ and $\psi_{i}$ only have nondegenerate critical points given
by $\left\{  x\in\Omega;\left\vert \nabla\psi_{i}\left(  x\right)  \right\vert
=0\right\}  =\left\{  p_{j};j=1,\cdot\cdot,d\right\}  $. Finally,
$\underset{\overline{\Omega}}{\text{max}}\psi_{i}=\underset{\overline{\Omega}%
}{\text{max}}\psi$ and $\psi\left(  a_{1}\right)  =\psi\left(  \left(
\Lambda_{i}\left(  1,\cdot\right)  \right)  ^{-1}\left(  \tau^{i-1}\left(
p_{1}\right)  \right)  \right)  =\psi\left(  \left(  \Lambda_{i}\left(
1,\cdot\right)  \right)  ^{-1}\left(  p_{i}\right)  \right)  =\psi_{i}\left(
p_{i}\right)  $ allow to conclude that $p_{i}$ is the unique global maximum of
$\psi_{i}$ and $\underset{\overline{\Omega}}{\text{max}}\psi_{j}%
=\underset{\overline{\Omega}}{\text{max}}\psi_{i}$ $\forall i,j$. This
completes the proof.

\bigskip

Our next result resume some identities linked to the Carleman commutator (see
\cite{P} and references therein).

\bigskip

Proposition 2 .- \textit{Let }%
\[
\Phi\left(  x,t\right)  =\frac{s\varphi\left(  x\right)  }{\Gamma\left(
t\right)  }\text{, }s>0\text{, }\Gamma\left(  t\right)  =T-t+h\text{,
}h>0\text{ \textit{and} }\varphi\in C^{\infty}\left(  \overline{\Omega
}\right)  \text{ .}%
\]
\textit{Define for any }$f\in H^{2}\left(  \Omega\right)  $\textit{ }%
\[
\left\{
\begin{array}
[c]{ll}%
\mathcal{A}_{\varphi}f=-\nabla\Phi\cdot\nabla f-\frac{1}{2}\Delta\Phi f\text{
,} & \\
\mathcal{S}_{\varphi}f=-\Delta f-\eta f\text{ \textit{where} }\eta=\frac{1}%
{2}\partial_{t}\Phi+\frac{1}{4}\left\vert \nabla\Phi\right\vert ^{2}\text{ ,}
& \\
\mathcal{S}_{\varphi}^{\prime}f=-\partial_{t}\eta f\text{ .} &
\end{array}
\right.
\]
\textit{Then we have}

\begin{description}
\item[$\left(  i\right)  $]
\[
\displaystyle\int_{\Omega}\mathcal{A}_{\varphi}ff=-\displaystyle\frac{1}%
{2}\int_{\partial\Omega}\partial_{n}\Phi\left\vert f\right\vert ^{2}%
\]

\item[$\left(  ii\right)  $]
\[
\displaystyle\int_{\Omega}\mathcal{S}_{\varphi}ff=\displaystyle\int_{\Omega
}\left\vert \nabla f\right\vert ^{2}-\displaystyle\int_{\Omega}\eta\left\vert
f\right\vert ^{2}-\displaystyle\int_{\partial\Omega}\partial_{n}ff
\]

\item[$\left(  iii\right)  $]
\[%
\begin{array}
[c]{ll}%
\displaystyle\int_{\Omega}\mathcal{S}_{\varphi}^{\prime}ff+2\displaystyle\int
_{\Omega}\mathcal{S}_{\varphi}f\mathcal{A}_{\varphi}f & =-2\displaystyle\int
_{\Omega}\nabla f\nabla^{2}\Phi\nabla f-\displaystyle\int_{\Omega}\nabla
f\Delta\nabla\Phi f\\
& \quad-\displaystyle\frac{2}{\Gamma}\int_{\Omega}\left(  \eta+\frac{1}%
{4}\left\vert \nabla\Phi\right\vert ^{2}+\frac{s}{4}\nabla\Phi\nabla
^{2}\varphi\nabla\Phi\right)  \left\vert f\right\vert ^{2}\\
& \quad+\text{Boundary terms}%
\end{array}
\]
\textit{where}
\[%
\begin{array}
[c]{ll}%
\text{Boundary terms} & =\displaystyle2\int_{\partial\Omega}\partial
_{n}f\nabla\Phi\cdot\nabla f-\displaystyle\int_{\partial\Omega}\partial
_{n}\Phi\left\vert \nabla f\right\vert ^{2}\\
& \quad+\displaystyle\int_{\partial\Omega}\partial_{n}f\Delta\Phi
f+\displaystyle\int_{\partial\Omega}\eta\partial_{n}\Phi\left\vert
f\right\vert ^{2}\text{ .}%
\end{array}
\]

\end{description}

\bigskip

Proof .- The proof of $\int_{\Omega}\mathcal{A}_{\varphi}ff$ and $\int
_{\Omega}\mathcal{S}_{\varphi}ff$ is quite clear by integrations by parts. Now
we compute the bracket $2\left\langle \mathcal{S}_{\varphi}f,\mathcal{A}%
_{\varphi}f\right\rangle $: We have from the definition of $\mathcal{S}%
_{\varphi}f$ and $\mathcal{A}_{\varphi}f$,
\[
2\left\langle \mathcal{S}_{\varphi}f,\mathcal{A}_{\varphi}f\right\rangle
=2\int_{\Omega}\left(  \Delta f+\eta f\right)  \left(  \nabla\Phi\cdot\nabla
f+\frac{1}{2}\Delta\Phi f\right)
\]
and four integrations by parts give%
\[
2\left\langle \mathcal{S}_{\varphi}f,\mathcal{A}_{\varphi}f\right\rangle
=-2\int_{\Omega}\nabla f\nabla^{2}\Phi\nabla f-\int_{\Omega}\nabla
f\Delta\nabla\Phi f-\int_{\Omega}\nabla\eta\cdot\nabla\Phi\left\vert
f\right\vert ^{2}+\text{Boundary terms .}%
\]
Indeed,
\[
\int_{\Omega}\Delta f\nabla\Phi\cdot\nabla f=\int_{\partial\Omega}\partial
_{n}f\nabla\Phi\cdot\nabla f-\int_{\Omega}\nabla f\nabla^{2}\Phi\nabla
f-\int_{\Omega}\nabla f\nabla^{2}f\nabla\Phi\text{ ,}%
\]
but
\[
\int_{\Omega}\nabla f\nabla^{2}f\nabla\Phi=\frac{1}{2}\int_{\partial\Omega
}\partial_{n}\Phi\left\vert \nabla f\right\vert ^{2}-\frac{1}{2}\int_{\Omega
}\Delta\Phi\left\vert \nabla f\right\vert ^{2}\text{ .}%
\]
Second,
\[
\int_{\Omega}\Delta f\Delta\Phi f=\int_{\partial\Omega}\partial_{n}f\Delta\Phi
f-\int_{\Omega}\nabla f\Delta\nabla\Phi f-\int_{\Omega}\Delta\Phi\left\vert
\nabla f\right\vert ^{2}\text{ .}%
\]
Third,
\[
2\int_{\Omega}\eta f\nabla\Phi\cdot\nabla f=\int_{\partial\Omega}\eta
\partial_{n}\Phi\left\vert f\right\vert ^{2}-\int_{\Omega}\nabla\eta
\cdot\nabla\Phi\left\vert f\right\vert ^{2}-\int_{\Omega}\eta\Delta
\Phi\left\vert f\right\vert ^{2}\text{ .}%
\]
This concludes to the identity%
\[%
\begin{array}
[c]{ll}%
2\displaystyle\int_{\Omega}\mathcal{S}_{\varphi}f\mathcal{A}_{\varphi
}f-\displaystyle\int_{\Omega}\partial_{t}\eta\left\vert f\right\vert ^{2} &
=-2\displaystyle\int_{\Omega}\nabla f\nabla^{2}\Phi\nabla f-\displaystyle\int
_{\Omega}\nabla f\Delta\nabla\Phi f\\
& \quad+\text{Boundary terms}+\displaystyle\int_{\Omega}\left(  -\partial
_{t}\eta-\nabla\eta\cdot\nabla\Phi\right)  \left\vert f\right\vert ^{2}\text{
.}%
\end{array}
\]
Finally, using $\partial_{t}\Phi=\frac{1}{\Gamma}\Phi$ and $\partial_{t}%
^{2}\Phi=\frac{2}{\Gamma}\partial_{t}\Phi$, we obtain
\[%
\begin{array}
[c]{ll}%
-\partial_{t}\eta-\nabla\eta\cdot\nabla\Phi & =-\frac{1}{2}\partial_{t}%
^{2}\Phi-\nabla\Phi\cdot\nabla\partial_{t}\Phi-\frac{1}{2}\nabla\Phi\nabla
^{2}\Phi\nabla\Phi\\
& =-\frac{1}{\Gamma}\partial_{t}\Phi-\frac{1}{\Gamma}\left\vert \nabla
\Phi\right\vert ^{2}-\frac{s}{2\Gamma}\nabla\Phi\nabla^{2}\varphi\nabla\Phi\\
& =-\frac{2}{\Gamma}\eta-\frac{1}{2\Gamma}\left\vert \nabla\Phi\right\vert
^{2}-\frac{s}{2\Gamma}\nabla\Phi\nabla^{2}\varphi\nabla\Phi\text{ .}%
\end{array}
\]
This completes the proof of $\left(  iii\right)  $.

\bigskip

Recall the following result which is a variant of \cite[Lemma 4.3]{BP}.

\bigskip

Proposition 3 .- \textit{Let }$h>0$\textit{, }$T>0$\textit{\ and }$F_{1}%
,F_{2}\geq0$\textit{. Consider two positive functions }$y,N\in C^{1}\left(
\left[  0,T\right]  \right)  $\textit{\ such that }%
\begin{equation}
\left\{
\begin{array}
[c]{ll}%
\left\vert \displaystyle\frac{1}{2}y^{\prime}\left(  t\right)  +N\left(
t\right)  y\left(  t\right)  \right\vert \leq F_{1}y\left(  t\right)  \text{
,} & \\
N^{\prime}\left(  t\right)  \leq\displaystyle\frac{1+C_{0}}{T-t+h}N\left(
t\right)  +F_{2}\text{ ,} &
\end{array}
\right.  \tag{2.1}\label{2.1}%
\end{equation}
\textit{where }$C_{0}\geq0$\textit{. Then for any }$0\leq t_{1}<t_{2}%
<t_{3}\leq T$\textit{, one has }%
\[
y\left(  t_{2}\right)  ^{1+M}\leq y\left(  t_{3}\right)  y\left(
t_{1}\right)  ^{M}e^{D}%
\]
\textit{with}%
\[
M=\frac{\displaystyle\int_{t_{2}}^{t_{3}}\frac{1}{\left(  T-t+h\right)
^{1+C_{0}}}dt}{\displaystyle\int_{t_{1}}^{t_{2}}\frac{1}{\left(  T-t+h\right)
^{1+C_{0}}}dt}%
\]
\textit{and}%
\[
D=2M\left(  F_{2}\left(  t_{3}-t_{1}\right)  ^{2}+F_{1}\left(  t_{3}%
-t_{1}\right)  \right)  \text{ .}%
\]

\bigskip

Proof .- Set $\Gamma\left(  t\right)  =T-t+h$. From the second inequality of
(\ref{2.1}), we have
\begin{equation}
\left(  \Gamma^{1+C_{0}}N\right)  ^{\prime}\leq F_{2}\Gamma^{1+C_{0}}\text{ .}
\tag{2.2}\label{2.2}%
\end{equation}
Integrating (\ref{2.2}) over $(t,t_{2})$ with $t\in(t_{1},t_{2})$ gives
\[
\left(  \frac{\Gamma\left(  t_{2}\right)  }{\Gamma\left(  t\right)  }\right)
^{1+C_{0}}N\left(  t_{2}\right)  \leq N\left(  t\right)  +F_{2}\left(
t_{2}-t_{1}\right)  \text{ .}%
\]
By the first inequality of (\ref{2.1}),
\[
y^{\prime}(t)+2N(t)y(t)\leq2F_{1}y(t)
\]
and we derive that
\[
y^{\prime}+\left(  2\left(  \frac{\Gamma\left(  t_{2}\right)  }{\Gamma\left(
t\right)  }\right)  ^{1+C_{0}}N\left(  t_{2}\right)  -2F_{2}\left(
t_{2}-t_{1}\right)  -2F_{1}\right)  y\leq0\text{ for }t\in(t_{1},t_{2})\text{
.}%
\]
Integrating over $\left(  t_{1},t_{2}\right)  $, we obtain
\begin{equation}
y(t_{2})e^{2N\left(  t_{2}\right)  \displaystyle\int_{t_{1}}^{t_{2}}\left(
\frac{\Gamma\left(  t_{2}\right)  }{\Gamma\left(  t\right)  }\right)
^{1+C_{0}}dt}\leq y\left(  t_{1}\right)  e^{2F_{2}\left(  t_{2}-t_{1}\right)
^{2}+2F_{1}\left(  t_{2}-t_{1}\right)  }\text{ .} \tag{2.3}\label{2.3}%
\end{equation}
On the other hand, integrating (\ref{2.2}) over $\left(  t_{2},t\right)  $
with $t\in(t_{2},t_{3})$, one has
\[
N\left(  t\right)  \leq\left(  \frac{\Gamma\left(  t_{2}\right)  }%
{\Gamma\left(  t\right)  }\right)  ^{1+C_{0}}\left(  N\left(  t_{2}\right)
+F_{2}\left(  t_{3}-t_{2}\right)  \right)  \text{ .}%
\]
By the first inequality of (\ref{2.1}),
\[
-y^{\prime}(t)-2N(t)y(t)\leq2F_{1}y(t)
\]
and it follows that
\[
0\leq y^{\prime}+\left[  2\left(  \frac{\Gamma\left(  t_{2}\right)  }%
{\Gamma\left(  t\right)  }\right)  ^{1+C_{0}}\left(  N\left(  t_{2}\right)
+F_{2}\left(  t_{3}-t_{2}\right)  \right)  +2F_{1}\right]  y\text{ for }%
t\in(t_{2},t_{3})\text{ .}%
\]
Integrating over $\left(  t_{2},t_{3}\right)  $ yields
\begin{equation}
y\left(  t_{2}\right)  \leq e^{2\left(  N\left(  t_{2}\right)  +F_{2}\left(
t_{3}-t_{2}\right)  \right)  \displaystyle\int_{t_{2}}^{t_{3}}\left(
\frac{\Gamma\left(  t_{2}\right)  }{\Gamma\left(  t\right)  }\right)
^{1+C_{0}}dt}y\left(  t_{3}\right)  e^{2F_{1}\left(  t_{3}-t_{2}\right)
}\text{ .} \tag{2.4}\label{2.4}%
\end{equation}
Combining (\ref{2.3}) and (\ref{2.4}), one has%
\[
y\left(  t_{2}\right)  \leq y\left(  t_{3}\right)  \left(  \frac{y\left(
t_{1}\right)  }{y\left(  t_{2}\right)  }e^{2F_{2}\left(  t_{2}-t_{1}\right)
^{2}}e^{2F_{1}\left(  t_{2}-t_{1}\right)  }\right)  ^{M}e^{2F_{1}\left(
t_{3}-t_{2}\right)  }e^{2F_{2}\left(  t_{3}-t_{2}\right)  \displaystyle\int
_{t_{2}}^{t_{3}}\left(  \frac{\Gamma\left(  t_{2}\right)  }{\Gamma\left(
t\right)  }\right)  ^{1+C_{0}}dt}%
\]
which gives%
\[
y\left(  t_{2}\right)  \leq y\left(  t_{3}\right)  \left(  \frac{y\left(
t_{1}\right)  }{y\left(  t_{2}\right)  }\right)  ^{M}e^{2F_{2}\left(
t_{2}-t_{1}\right)  ^{2}M}e^{2F_{1}\left(  t_{2}-t_{1}\right)  M}%
e^{2F_{1}\left(  t_{3}-t_{2}\right)  }e^{2F_{2}\left(  t_{3}-t_{2}\right)
\left(  t_{2}-t_{1}\right)  M}%
\]
which implies the desired estimate since $M>1$.

\bigskip

\section{Proof of Theorem 1}

\bigskip

The plan of the proof of Theorem 1 is as follows. We divide it into seven
steps. In Step 1, we derive some estimates on the Morse functions given in
Proposition 1. In Step 2, we introduce the weight functions and establish the
key properties linked to the Morse functions. In Step 3, we perform a change
of function and introduce the operators described in Proposition 2. In Step 4,
we construct a new frequency function adapted to our global approach. In Step
5, key estimates for the Carleman operator is provided. In Step 6, we solve a
system of ordinary differential inequalities thanks to Proposition 3. In Step
7, we conclude the proof by making appear the control domain $\omega
\times\left\{  T\right\}  $.

\bigskip

\subsection{Step 1: The Morse functions}

\bigskip

We have by Proposition 1, the existence of Morse functions $\psi_{i}$
associated to a critical point $p_{i}$ which is its unique global maximum in
$\overline{\Omega}$. By Morse Lemma, there exists a neighborhood of $p_{i}$
and a diffeomorphism $U$ such that $U\left(  p_{i}\right)  =0$ and locally
\[
\psi_{i}\left(  U^{-1}\left(  x\right)  \right)  =\psi_{i}\left(
p_{i}\right)  -\left\vert x\right\vert ^{2}%
\]
which implies
\[
\frac{1}{4}\left\vert \text{Jac}U^{-1}\left(  x\right)  \nabla\psi_{i}\left(
U^{-1}\left(  x\right)  \right)  \right\vert ^{2}=\left\vert x\right\vert
^{2}=\psi_{i}\left(  p_{i}\right)  -\psi_{i}\left(  U^{-1}\left(  x\right)
\right)
\]
and consequently, there are $c_{1},c_{2}>0$ such that for any $i\in\left\{
1,\cdot\cdot,d\right\}  $, in a neighborhood of $p_{i}$
\begin{equation}
c_{1}\left\vert \nabla\psi_{i}\right\vert ^{2}\leq\left(  \underset
{\overline{\Omega}}{\text{max}}\psi_{i}-\psi_{i}\right)  \leq c_{2}\left\vert
\nabla\psi_{i}\right\vert ^{2}\text{ .} \tag{3.1.1}\label{3.1.1}%
\end{equation}

\bigskip

Let
\[
\mathcal{B}_{i}\text{ be a neighborhood of }\left\{  x\in\Omega;\left\vert
\nabla\psi_{i}\left(  x\right)  \right\vert =0\text{ and }\underset
{\overline{\Omega}}{\text{max}}\psi_{i}-\psi_{i}\left(  x\right)  =0\right\}
\]
in which (\ref{3.1.1}) holds ,%
\[
\mathcal{C}_{i}\text{ be a neighborhood of }\left\{  x\in\Omega;\left\vert
\nabla\psi_{i}\left(  x\right)  \right\vert =0\text{ and }\underset
{\overline{\Omega}}{\text{max}}\psi_{i}-\psi_{i}\left(  x\right)  >0\right\}
\]
with $\mathcal{B}_{i}\cap\mathcal{C}_{i}=\emptyset$ in which $\psi_{i}%
-\psi_{j}<0$ for some $j\neq i$. This is possible because $\psi_{i}\left(
p_{j}\right)  <\psi_{j}\left(  p_{j}\right)  $ using Proposition 1 $\left(
iv\right)  $ and $\left(  v\right)  $ with $\left\{  p_{j};j=1,\cdot
\cdot,d\right\}  =\left\{  x\in\Omega;\left\vert \nabla\psi_{i}\left(
x\right)  \right\vert =0\right\}  $ and $\mathcal{C}_{i}=%
{\textstyle\bigcup\limits_{j\neq i}}
\Theta_{p_{j}}$ where $\Theta_{p_{j}}$ is a sufficiently small neighborhood of
$p_{j}$. And finally let%
\[
\mathcal{D}_{i}=\Omega\left\backslash \left(  \mathcal{B}_{i}\cup
\mathcal{C}_{i}\right)  \right.  \text{ be such that }\Omega=\mathcal{B}%
_{i}\cup\mathcal{C}_{i}\cup\mathcal{D}_{i}\text{ .}%
\]

\bigskip

Proposition 4 .- \textit{There are }$c_{1}>0$\textit{ and }$c_{2}>0$\textit{
such that for any }$i\in\left\{  1,\cdot\cdot,d\right\}  $

\begin{description}
\item[$\left(  i\right)  $] \textit{In }$D_{i}$\textit{, }%
\[
c_{1}\left\vert \nabla\psi_{i}\right\vert ^{2}\leq\left(  \underset
{\overline{\Omega}}{\text{max}}\psi_{i}-\psi_{i}\right)  \leq c_{2}\left\vert
\nabla\psi_{i}\right\vert ^{2}\text{ .}%
\]

\item[$\left(  ii\right)  $] \textit{In }$B_{i}$\textit{,}%
\[
c_{1}\left\vert \nabla\psi_{i}\right\vert ^{2}\leq\left(  \underset
{\overline{\Omega}}{\text{max}}\psi_{i}-\psi_{i}\right)  \leq c_{2}\left\vert
\nabla\psi_{i}\right\vert ^{2}\text{ .}%
\]

\item[$\left(  iii\right)  $] \textit{In }$C_{i}$\textit{,}%
\[
c_{1}\left\vert \nabla\psi_{i}\right\vert ^{2}\leq\left(  \underset
{\overline{\Omega}}{\text{max}}\psi_{i}-\psi_{i}\right)  \text{ .}%
\]

\end{description}

\bigskip

Proof .- The inequality $\left(  ii\right)  $ holds by definition of
$\mathcal{B}_{i}$. In $\mathcal{C}_{i}$, we use $\underset{\overline{\Omega}%
}{\text{max}}\psi_{i}-\psi_{i}\geq c>0$ and $\left\vert \nabla\psi
_{i}\right\vert ^{2}\leq\underset{\overline{\Omega}}{\text{max}}\left\vert
\nabla\psi_{i}\right\vert ^{2}\leq\frac{\underset{\overline{\Omega}%
}{\text{max}}\left\vert \nabla\psi_{i}\right\vert ^{2}}{c}\left(
\underset{\overline{\Omega}}{\text{max}}\psi_{i}-\psi_{i}\right)  $. In
$\mathcal{D}_{i}$, $\left\vert \nabla\psi_{i}\right\vert >0$ and
$\underset{\overline{\Omega}}{\text{max}}\psi_{i}-\psi_{i}>0\ $imply the
desired estimates.

\bigskip

\subsection{Step 2: The weight functions}

\bigskip

Introduce for any $i\in\left\{  1,\cdot\cdot,d\right\}  $%
\[
\left\{
\begin{array}
[c]{ll}%
\varphi_{i,1}=\psi_{i}-\underset{\overline{\Omega}}{\text{max}}\psi_{i}\text{
,} & \\
\varphi_{i,2}=-\psi_{i}-\underset{\overline{\Omega}}{\text{max}}\psi_{i}\text{
.} &
\end{array}
\right.
\]
Notice that
\begin{equation}
\varphi_{i,1}=\varphi_{i,2}\text{ on }\partial\Omega\text{ and }\partial
_{n}\varphi_{i,1}+\partial_{n}\varphi_{i,2}=0\text{ on }\partial\Omega\text{
.} \tag{3.2.1}\label{3.2.1}%
\end{equation}
Further, the link between $\varphi_{i,1}$ and $\psi_{i}$ is described as
follows: $\left\vert \varphi_{i,1}\right\vert =\underset{\overline{\Omega}%
}{\text{max}}\psi_{i}-\psi_{i}$ and $\left\vert \nabla\varphi_{i,1}\right\vert
^{2}=\left\vert \nabla\psi_{i}\right\vert ^{2}$. Now, we are able to state the
properties of $\varphi_{i,1}$ and $\varphi_{i,2}$.

\bigskip

Proposition 5 .- \textit{There are }$c_{1},\cdot\cdot,c_{6}>0$\textit{ all
positive constants such that for any }$i\in\left\{  1,\cdot\cdot,d\right\}  $

\begin{description}
\item[$\left(  i\right)  $] \textit{In }$D_{i}$\textit{,}%
\[
c_{1}\left\vert \nabla\varphi_{i,1}\right\vert ^{2}\leq\left\vert
\varphi_{i,1}\right\vert \leq c_{2}\left\vert \nabla\varphi_{i,1}\right\vert
^{2}\text{ .}%
\]

\item[$\left(  ii\right)  $] \textit{In }$B_{i}$\textit{,}%
\[
c_{1}\left\vert \nabla\varphi_{i,1}\right\vert ^{2}\leq\left\vert
\varphi_{i,1}\right\vert \leq c_{2}\left\vert \nabla\varphi_{i,1}\right\vert
^{2}\text{ .}%
\]

\item[$\left(  iii\right)  $] \textit{In }$C_{i}$\textit{,}%
\[
c_{1}\left\vert \nabla\varphi_{i,1}\right\vert ^{2}\leq\left\vert
\varphi_{i,1}\right\vert \text{ .}%
\]

\item[$\left(  iv\right)  $] \textit{There is }$j\in\left\{  1,\cdot
\cdot,d\right\}  $\textit{ with }$j\neq i$\textit{ such that}%
\[
\varphi_{i,1}-\varphi_{j,1}\leq-c_{3}\text{ \textit{in} }\mathcal{C}_{i}\text{
.}%
\]

\item[$\left(  v\right)  $]
\[
c_{4}\left\vert \nabla\varphi_{i,2}\right\vert ^{2}\leq\left\vert
\varphi_{i,2}\right\vert \text{ \textit{in} }\Omega\text{ \textit{and}
}\left\vert \varphi_{i,2}\right\vert \leq c_{5}\left\vert \nabla\varphi
_{i,2}\right\vert ^{2}\text{ \textit{in a neighborhood of} }\partial
\Omega\text{ .}%
\]

\item[$\left(  vi\right)  $]
\[
\varphi_{i,2}-\varphi_{i,1}\leq-c_{6}\text{ \textit{outside a neighborhood of}
}\partial\Omega\text{ .}%
\]

\end{description}

\bigskip

Proof .- By the properties of the Morse functions described in Proposition 4,
we deduce $\left(  i\right)  -\left(  ii\right)  $ and $\left(  iii\right)  $.
The inequality $\left(  iv\right)  $ holds from the definition of
$\mathcal{C}_{i}$ and Proposition 1 $\left(  v\right)  $. Next, we start to
prove $\left(  v\right)  $ by seeing that $\left\vert \nabla\varphi
_{i,2}\right\vert ^{2}\leq c\leq\frac{c}{\underset{\overline{\Omega}%
}{\text{max}}\psi_{i}}\left\vert \varphi_{i,2}\right\vert $. Since $\left\vert
\nabla\varphi_{i,2}\right\vert =\left\vert \nabla\psi_{i}\right\vert >0$ in a
neighborhood of $\partial\Omega$, we have $\left\vert \varphi_{i,2}\right\vert
\leq c\leq c_{5}\left\vert \nabla\varphi_{i,2}\right\vert ^{2}$. This
completes the proof of $\left(  v\right)  $. Finally, since $\psi_{i}>0$
outside a neighborhood of $\partial\Omega$, we get $0<c\leq\psi_{i}$ and
$\varphi_{i,2}-\varphi_{i,1}=-2\psi_{i}\leq-2c=-c_{6}$, that is $\left(
vi\right)  $.

\bigskip

\subsection{Step 3: Change of functions}

\bigskip

Introduce for any $\left(  x,t\right)  \in\Omega\times\left[  0,T\right]  $
and any $i\in\left\{  1,\cdot\cdot,d\right\}  $
\[
\left\{
\begin{array}
[c]{ll}%
\Phi_{i}\left(  x,t\right)  =\displaystyle\frac{s}{\Gamma\left(  t\right)
}\varphi_{i,1}\left(  x\right)  \text{ ,} & \\
\Phi_{d+i}\left(  x,t\right)  =\displaystyle\frac{s}{\Gamma\left(  t\right)
}\varphi_{i,2}\left(  x\right)  \text{ .} &
\end{array}
\right.
\]
with $s\in\left(  0,1\right]  $ and $\Gamma\left(  t\right)  =T-t+h$,
$h\in\left(  0,1\right]  $.

Let $\boldsymbol{f}=\left(  f_{i}\right)  _{1\leq i\leq2d}$ where
$f_{i}=ue^{\Phi_{i}/2}$. We look for the equation solved by $f_{i}$ by
computing $e^{\Phi_{i}/2}\left(  \partial_{t}-\Delta\right)  \left(
e^{-\Phi_{i}/2}f_{i}\right)  $. Introduce
\[
\left\{
\begin{array}
[c]{ll}%
\mathcal{A}_{\varphi_{i}}f_{i}=-\nabla\Phi_{i}\cdot\nabla f_{i}-\frac{1}%
{2}\Delta\Phi_{i}f_{i}\text{ ,} & \\
\mathcal{S}_{\varphi_{i}}f_{i}=-\Delta f_{i}-\eta_{i}f_{i}\text{ where }%
\eta_{i}=\frac{1}{2}\partial_{t}\Phi_{i}+\frac{1}{4}\left\vert \nabla\Phi
_{i}\right\vert ^{2}\text{ .} &
\end{array}
\right.
\]
Let $\mathcal{S}\boldsymbol{f}=\left(  \mathcal{S}_{\varphi_{i}}f_{i}\right)
_{1\leq i\leq2d}$, $\mathcal{A}\boldsymbol{f}=\left(  \mathcal{A}_{\varphi
_{i}}f_{i}\right)  _{1\leq i\leq2d}$, and $\digamma=\left(  -af_{i}\right)
_{1\leq i\leq2d}$. We find that
\begin{equation}
\left\{
\begin{array}
[c]{ll}%
\partial_{t}\boldsymbol{f}+\mathcal{S}\boldsymbol{f}=\mathcal{A}%
\boldsymbol{f}+\digamma\text{ ,} & \\
\partial_{n}f_{i}-\frac{1}{2}\partial_{n}\Phi_{i}f_{i}=0\text{ on }%
\partial\Omega\times\left(  0,T\right)  \text{ .} &
\end{array}
\right.  \tag{3.3.1}\label{3.3.1}%
\end{equation}
Let $\left\langle \cdot,\cdot\right\rangle $ denote the usual scalar product
in $\left(  L^{2}\left(  \Omega\right)  \right)  ^{2d}$ and let $\left\Vert
\cdot\right\Vert $ be its corresponding norm. Now, we claim that
\begin{equation}
\left\{
\begin{array}
[c]{ll}%
\left\langle \mathcal{A}\boldsymbol{f},\boldsymbol{f}\right\rangle =0\text{ ,}
& \\
\left\langle \mathcal{S}\boldsymbol{f},\boldsymbol{f}\right\rangle
=\displaystyle\sum_{i=1,..,2d}\displaystyle\int_{\Omega}\left\vert \nabla
f_{i}\right\vert ^{2}-\displaystyle\int_{\Omega}\eta_{i}\left\vert
f_{i}\right\vert ^{2}\text{ ,} & \\
\displaystyle\frac{d}{dt}\left\langle \mathcal{S}\boldsymbol{f},\boldsymbol{f}%
\right\rangle =-\displaystyle\sum_{i=1,..,2d}\displaystyle\int_{\Omega
}\partial_{t}\eta_{i}\left\vert f_{i}\right\vert ^{2}+2\left\langle
\mathcal{S}\boldsymbol{f},\partial_{t}\boldsymbol{f}\right\rangle
:=\left\langle \mathcal{S}^{\prime}\boldsymbol{f},\boldsymbol{f}\right\rangle
+2\left\langle \mathcal{S}\boldsymbol{f},\partial_{t}\boldsymbol{f}%
\right\rangle \text{ .} &
\end{array}
\right.  \tag{3.3.2}\label{3.3.2}%
\end{equation}
Indeed, applying Proposition 2 $\left(  i\right)  -\left(  ii\right)  $ and
using the Robin boundary condition for $f_{i}$, all the boundary terms
appearing in the integrations by parts can be dropped since for any
$i\in\left\{  1,\cdot\cdot,d\right\}  $
\begin{equation}
\Phi_{i}=\Phi_{d+i}\text{ and }\partial_{n}\Phi_{i}+\partial_{n}\Phi
_{d+i}=0\text{ on }\partial\Omega\times\left(  0,T\right)  \text{ ,}
\tag{3.3.3}\label{3.3.3}%
\end{equation}
by (\ref{3.2.1}). To establish the last identity in (\ref{3.3.2}), we compute
$\frac{d}{dt}\left\langle \mathcal{S}\boldsymbol{f},\boldsymbol{f}%
\right\rangle $ as follows:%
\[%
\begin{array}
[c]{ll}%
\displaystyle\frac{d}{dt}\left\langle \mathcal{S}\boldsymbol{f},\boldsymbol{f}%
\right\rangle  & =\displaystyle\frac{d}{dt}\left(  \displaystyle\sum
_{i=1,..,2d}\displaystyle\int_{\Omega}\left\vert \nabla f_{i}\right\vert
^{2}-\displaystyle\int_{\Omega}\eta_{i}\left\vert f_{i}\right\vert ^{2}\right)
\\
& =2\left\langle \mathcal{S}\boldsymbol{f},\partial_{t}\boldsymbol{f}%
\right\rangle -\displaystyle\sum_{i=1,..,2d}\displaystyle\int_{\Omega}%
\partial_{t}\eta_{i}\left\vert f_{i}\right\vert ^{2}+\displaystyle2\sum
_{i=1,..,2d}\displaystyle\int_{\partial\Omega}\partial_{n}f_{i}\partial
_{t}f_{i}%
\end{array}
\]
by an integration by parts. But, by using the Robin boundary condition for
$f_{i}=ue^{\Phi_{i}/2}$ in (\ref{3.3.1}), we have%
\[
\sum_{i=1,..,2d}\int_{\partial\Omega}\partial_{n}f_{i}\partial_{t}f_{i}%
=\sum_{i=1,..,2d}\int_{\partial\Omega}\frac{1}{2}\partial_{n}\Phi_{i}\left(
u\partial_{t}u+\left\vert u\right\vert ^{2}\frac{1}{2}\partial_{t}\Phi
_{i}\right)  e^{\Phi_{i}}=0
\]
since for any $i\in\left\{  1,\cdot\cdot,d\right\}  $, $\Phi_{d+i}=\Phi_{i}$
and $\partial_{n}\Phi_{i}+\partial_{n}\Phi_{d+i}=0$ on $\partial\Omega
\times\left(  0,T\right)  $.

\bigskip

\subsection{Step 4: Energy estimates}

\bigskip

By a standard energy method, we have%
\[
\frac{1}{2}\frac{d}{dt}\left\Vert \boldsymbol{f}\right\Vert ^{2}+\left\langle
\mathcal{S}\boldsymbol{f},\boldsymbol{f}\right\rangle =\left\langle
\digamma,\boldsymbol{f}\right\rangle \text{ ,}%
\]
and by introducing the frequency function
\[
\mathbf{N}\left(  t\right)  =\frac{\left\langle \mathcal{S}\boldsymbol{f}%
,\boldsymbol{f}\right\rangle }{\left\Vert \boldsymbol{f}\right\Vert ^{2}}%
\]
it holds%
\[
\mathbf{N}^{\prime}\left(  t\right)  \left\Vert \boldsymbol{f}\right\Vert
^{2}\leq\left\langle \mathcal{S}^{\prime}\boldsymbol{f},\boldsymbol{f}%
\right\rangle +2\left\langle \mathcal{S}\boldsymbol{f},\mathcal{A}%
\boldsymbol{f}\right\rangle +\left\Vert \digamma\right\Vert ^{2}\text{ .}%
\]
Indeed, for the energy identity we use the first equality of (\ref{3.3.1}) and
$\left\langle \mathcal{A}\boldsymbol{f},\boldsymbol{f}\right\rangle =0$. For
the inequality of the derivative of the frequency function, we use $\frac
{d}{dt}\left\langle \mathcal{S}\boldsymbol{f},\boldsymbol{f}\right\rangle
=\left\langle \mathcal{S}^{\prime}\boldsymbol{f},\boldsymbol{f}\right\rangle
+2\left\langle \mathcal{S}\boldsymbol{f},\partial_{t}\boldsymbol{f}%
\right\rangle $ (see (\ref{3.3.2})) and replace $\partial_{t}\boldsymbol{f}$
by $\mathcal{A}\boldsymbol{f}-\mathcal{S}\boldsymbol{f}+\digamma$ in order to
get
\[%
\begin{array}
[c]{ll}%
\mathbf{N}^{\prime}\left(  t\right)  \left\Vert \boldsymbol{f}\right\Vert ^{4}
& =\left(  \left\langle \mathcal{S}^{\prime}\boldsymbol{f},\boldsymbol{f}%
\right\rangle +2\left\langle \mathcal{S}\boldsymbol{f},\partial_{t}%
\boldsymbol{f}\right\rangle \right)  \left\Vert \boldsymbol{f}\right\Vert
^{2}-\left\langle \mathcal{S}\boldsymbol{f},\boldsymbol{f}\right\rangle
\left(  -2\left\langle \mathcal{S}\boldsymbol{f},\boldsymbol{f}\right\rangle
+2\left\langle \digamma,\boldsymbol{f}\right\rangle \right) \\
& =\left(  \left\langle \mathcal{S}^{\prime}\boldsymbol{f},\boldsymbol{f}%
\right\rangle +2\left\langle \mathcal{S}\boldsymbol{f},\mathcal{A}%
\boldsymbol{f}\right\rangle \right)  \left\Vert \boldsymbol{f}\right\Vert
^{2}-2\left\Vert \mathcal{S}\boldsymbol{f}\right\Vert ^{2}\left\Vert
\boldsymbol{f}\right\Vert ^{2}+2\left\langle \mathcal{S}\boldsymbol{f}%
,\digamma\right\rangle \left\Vert \boldsymbol{f}\right\Vert ^{2}\\
& \quad+2\left\langle \mathcal{S}\boldsymbol{f},\boldsymbol{f}\right\rangle
^{2}-2\left\langle \mathcal{S}\boldsymbol{f},\boldsymbol{f}\right\rangle
\left\langle \digamma,\boldsymbol{f}\right\rangle \\
& =\left(  \left\langle \mathcal{S}^{\prime}\boldsymbol{f},\boldsymbol{f}%
\right\rangle +2\left\langle \mathcal{S}\boldsymbol{f},\mathcal{A}%
\boldsymbol{f}\right\rangle \right)  \left\Vert \boldsymbol{f}\right\Vert
^{2}-2\left\Vert \mathcal{S}\boldsymbol{f}-\frac{1}{2}\digamma\right\Vert
^{2}\left\Vert \boldsymbol{f}\right\Vert ^{2}+\frac{1}{2}\left\Vert
\digamma\right\Vert ^{2}\left\Vert \boldsymbol{f}\right\Vert ^{2}\\
& \quad+2\left\langle \mathcal{S}\boldsymbol{f}-\frac{1}{2}\digamma
,\boldsymbol{f}\right\rangle ^{2}-\frac{1}{2}\left\langle \digamma
,\boldsymbol{f}\right\rangle ^{2}\text{ .}%
\end{array}
\]
By Cauchy-Schwarz, we obtain the desired estimate for $\mathbf{N}^{\prime
}\left(  t\right)  $.

\bigskip

Since
\[
\left\Vert \digamma\right\Vert ^{2}\leq\left\Vert a\right\Vert _{\infty}%
^{2}\left\Vert \boldsymbol{f}\right\Vert ^{2}%
\]
where $\left\Vert a\right\Vert _{\infty}=\left\Vert a\right\Vert _{L^{\infty
}\left(  \Omega\times\left(  0,T\right)  \right)  }$, we obtain the following
system of ordinary differential inequalities%
\begin{equation}
\left\{
\begin{array}
[c]{ll}%
\left\vert \displaystyle\frac{1}{2}\frac{d}{dt}\left\Vert \boldsymbol{f}%
\right\Vert ^{2}+\mathbf{N}\left(  t\right)  \left\Vert \boldsymbol{f}%
\right\Vert ^{2}\right\vert \leq\left\Vert a\right\Vert _{\infty}\left\Vert
\boldsymbol{f}\right\Vert ^{2}\text{ ,} & \\
\mathbf{N}^{\prime}\left(  t\right)  \leq\displaystyle\frac{\left\langle
\mathcal{S}^{\prime}\boldsymbol{f},\boldsymbol{f}\right\rangle +2\left\langle
\mathcal{S}\boldsymbol{f},\mathcal{A}\boldsymbol{f}\right\rangle }{\left\Vert
\boldsymbol{f}\right\Vert ^{2}}+\left\Vert a\right\Vert _{\infty}^{2}\text{ .}
&
\end{array}
\right.  \tag{3.4.1}\label{3.4.1}%
\end{equation}

\bigskip

\subsection{Step 5: Carleman commutator estimates}

\bigskip

We claim that for some $s\in\left(  0,1\right]  $ sufficiently small,
$\eta_{i}\leq0$ and $\left\langle \mathcal{S}\boldsymbol{f},\boldsymbol{f}%
\right\rangle \geq0$ and
\[
\left\langle \mathcal{S}^{\prime}\boldsymbol{f},\boldsymbol{f}\right\rangle
+2\left\langle \mathcal{S}\boldsymbol{f},\mathcal{A}\boldsymbol{f}%
\right\rangle \leq\frac{1+C_{0}}{\Gamma}\left\langle \mathcal{S}%
\boldsymbol{f},\boldsymbol{f}\right\rangle +\frac{C}{h^{2}}\left\Vert
\boldsymbol{f}\right\Vert ^{2}\text{ ,}%
\]
where $C_{0}\in\left(  0,1\right)  $ and $C>0$ do not depend on $h\in\left(
0,1\right]  $.

\bigskip

Indeed, observe that
\[
\eta_{i}=\frac{1}{2}\partial_{t}\Phi_{i}+\frac{1}{4}\left\vert \nabla\Phi
_{i}\right\vert ^{2}=\left\vert
\begin{array}
[c]{ll}%
\frac{s}{4\Gamma^{2}}\left(  -2\left\vert \varphi_{i,1}\right\vert
+s\left\vert \nabla\varphi_{i,1}\right\vert ^{2}\right)  & \text{if }%
i\in\left\{  1,\cdot\cdot,d\right\} \\
\frac{s}{4\Gamma^{2}}\left(  -2\left\vert \varphi_{i-d,2}\right\vert
+s\left\vert \nabla\varphi_{i-d,2}\right\vert ^{2}\right)  & \text{if }%
i\in\left\{  d+1,\cdot\cdot,2d\right\}
\end{array}
\right.  \leq0
\]
for $s\in\left(  0,1\right]  $ sufficiently small since $\left\vert
\nabla\varphi_{i,j}\right\vert ^{2}\leq c\left\vert \varphi_{i,j}\right\vert $
for any $i\in\left\{  1,\cdot\cdot,d\right\}  $, any $j\in\left\{
1,2\right\}  $ by Proposition 5 $\left(  i\right)  -\left(  iii\right)  $ and
$\left(  v\right)  $. This concludes the proof that $\left\langle
\mathcal{S}\boldsymbol{f},\boldsymbol{f}\right\rangle \geq0$ for $s$ small.

\bigskip

By Proposition 2 $\left(  iii\right)  $,
\begin{equation}%
\begin{array}
[c]{ll}%
\left\langle \mathcal{S}^{\prime}\boldsymbol{f},\boldsymbol{f}\right\rangle
+2\left\langle \mathcal{S}\boldsymbol{f},\mathcal{A}\boldsymbol{f}%
\right\rangle  & =-2\displaystyle\sum_{i=1,..,2d}\displaystyle\int_{\Omega
}\nabla f_{i}\nabla^{2}\Phi_{i}\nabla f_{i}-\displaystyle\sum_{i=1,..,2d}%
\int_{\Omega}\nabla f_{i}\Delta\nabla\Phi_{i}f_{i}\\
& \quad-\displaystyle\frac{2}{\Gamma}\sum_{i=1,..,2d}\displaystyle\int
_{\Omega}\left(  \eta_{i}+\frac{1}{4}\left\vert \nabla\Phi_{i}\right\vert
^{2}+\frac{s}{4}\nabla\Phi_{i}\nabla^{2}\varphi_{i}\nabla\Phi_{i}\right)
\left\vert f_{i}\right\vert ^{2}\\
& \quad+\text{Boundary terms}%
\end{array}
\tag{3.5.1}\label{3.5.1}%
\end{equation}
where $\varphi_{i}=\varphi_{i,1}$ for $i\in\left\{  1,\cdot\cdot,d\right\}  $,
$\varphi_{i}=\varphi_{i-d,2}$ for $i\in\left\{  d+1,\cdot\cdot,2d\right\}  $,
and
\begin{equation}%
\begin{array}
[c]{ll}%
\text{Boundary terms} & =\displaystyle2\sum_{i=1,..,2d}\displaystyle\int
_{\partial\Omega}\partial_{n}f_{i}\nabla\Phi_{i}\cdot\nabla f_{i}%
-\displaystyle\sum_{i=1,..,2d}\displaystyle\int_{\partial\Omega}\partial
_{n}\Phi_{i}\left\vert \nabla f_{i}\right\vert ^{2}\\
& \quad+\displaystyle\sum_{i=1,..,2d}\displaystyle\int_{\partial\Omega
}\partial_{n}f_{i}\Delta\Phi_{i}f_{i}+\displaystyle\sum_{i=1,..,2d}%
\displaystyle\int_{\partial\Omega}\eta_{i}\partial_{n}\Phi_{i}\left\vert
f_{i}\right\vert ^{2}\text{ .}%
\end{array}
\tag{3.5.2}\label{3.5.2}%
\end{equation}

\bigskip

First we estimate the contribution of the gradient terms:%
\begin{equation}%
\begin{array}
[c]{ll}%
\displaystyle\sum_{i=1,..,2d}\left(  -2\displaystyle\int_{\Omega}\nabla
f_{i}\nabla^{2}\Phi_{i}\nabla f_{i}-\displaystyle\int_{\Omega}\nabla
f_{i}\Delta\nabla\Phi_{i}f_{i}\right)  & \leq\displaystyle\frac{cs}{\Gamma
}\sum_{i=1,..,2d}\int_{\Omega}\left\vert \nabla f_{i}\right\vert ^{2}%
+\frac{cs}{\Gamma}\left\Vert \boldsymbol{f}\right\Vert ^{2}\\
& \leq\displaystyle\frac{cs}{\Gamma}\sum_{i=1,..,2d}\int_{\Omega}\left\vert
\nabla f_{i}\right\vert ^{2}+\displaystyle\frac{c}{h}\left\Vert \boldsymbol{f}%
\right\Vert ^{2}\text{ }%
\end{array}
\tag{3.5.3}\label{3.5.3}%
\end{equation}
for $s\in\left(  0,1\right]  $, using Cauchy-Schwarz, $\left\vert 2\nabla
^{2}\Phi_{i}\right\vert \leq\frac{cs}{\Gamma}$, and $\left\vert \Delta
\nabla\Phi_{i}\right\vert \leq\frac{cs}{\Gamma}\leq\frac{1}{h}$.

Next we check the contribution of the boundary terms. We claim that
\[
\sum_{i=1,..,2d}\int_{\partial\Omega}\eta_{i}\partial_{n}\Phi_{i}\left\vert
f_{i}\right\vert ^{2}=0\text{ .}%
\]
Indeed, $\eta_{i}=\frac{1}{2}\partial_{t}\Phi_{i}+\frac{1}{4}\left\vert
\nabla\Phi_{i}\right\vert ^{2}$ implies
\[%
\begin{array}
[c]{ll}%
\displaystyle\sum_{i=1,..,2d}\displaystyle\int_{\partial\Omega}\eta
_{i}\partial_{n}\Phi_{i}\left\vert f_{i}\right\vert ^{2} & =\displaystyle\sum
_{i=1,..,d}\displaystyle\int_{\partial\Omega}\left(  \frac{1}{2}\partial
_{t}\Phi_{i}+\frac{1}{4}\left\vert \nabla\Phi_{i}\right\vert ^{2}\right)
\partial_{n}\Phi_{i}\left\vert u\right\vert ^{2}e^{\Phi_{i}}\\
& \quad+\displaystyle\sum_{i=1,..,d}\displaystyle\int_{\partial\Omega}\left(
\frac{1}{2}\partial_{t}\Phi_{i}+\frac{1}{4}\left\vert \nabla\Phi
_{i}\right\vert ^{2}\right)  \partial_{n}\Phi_{d+i}\left\vert u\right\vert
^{2}e^{\Phi_{i}}%
\end{array}
\]
where we used $\Phi_{d+i}=\Phi_{i}$ on $\partial\Omega\times\left(
0,T\right)  $ and $\left\vert \nabla\Phi_{d+i}\right\vert =\left\vert
\nabla\Phi_{i}\right\vert $ on $\partial\Omega\times\left(  0,T\right)  $.
Since $\partial_{n}\Phi_{i}+\partial_{n}\Phi_{d+i}=0$ on $\partial\Omega
\times\left(  0,T\right)  $, this completes the claim. We also have
\[
2\sum_{i=1,..,2d}\int_{\partial\Omega}\partial_{n}f_{i}\nabla\Phi_{i}%
\cdot\nabla f_{i}-\sum_{i=1,..,2d}\int_{\partial\Omega}\partial_{n}\Phi
_{i}\left\vert \nabla f_{i}\right\vert ^{2}=0\text{ .}%
\]
Indeed, since $\nabla\Phi_{i}=\partial_{n}\Phi_{i}\overrightarrow{n}$ on
$\partial\Omega\times\left(  0,T\right)  $ and $\partial_{n}f_{i}=\frac{1}%
{2}\partial_{n}\Phi_{i}f_{i}$,%
\[%
\begin{array}
[c]{ll}%
2\displaystyle\sum_{i=1,..,2d}\displaystyle\int_{\partial\Omega}\partial
_{n}f_{i}\nabla\Phi_{i}\cdot\nabla f_{i} & =2\displaystyle\sum_{i=1,..,2d}%
\displaystyle\int_{\partial\Omega}\partial_{n}\Phi_{i}\left\vert \frac{1}%
{2}\partial_{n}\Phi_{i}f_{i}\right\vert ^{2}\\
& =2\displaystyle\sum_{i=1,..,d}\displaystyle\int_{\partial\Omega}\left(
\partial_{n}\Phi_{i}+\partial_{n}\Phi_{d+i}\right)  \left\vert \frac{1}%
{2}\partial_{n}\Phi_{i}f_{i}\right\vert ^{2}\\
& =0
\end{array}
\]
where we used (\ref{3.3.3}). For the second contribution, it holds
\[
\left\vert \nabla f_{i}\right\vert ^{2}=\left\vert \nabla u+u\frac{1}{2}%
\nabla\Phi_{i}\right\vert ^{2}e^{\Phi_{i}}=\left\vert \partial_{\tau
}u\overrightarrow{\tau}+u\frac{1}{2}\partial_{n}\Phi_{i}\overrightarrow
{n}\right\vert ^{2}e^{\Phi_{i}}=\left(  \left\vert \partial_{\tau}u\right\vert
^{2}+\left\vert \frac{1}{2}u\partial_{n}\Phi_{i}\right\vert ^{2}\right)
e^{\Phi_{i}}%
\]
on $\partial\Omega\times\left(  0,T\right)  $. We then conclude that
$-\displaystyle\sum_{i=1,..,2d}\displaystyle\int_{\partial\Omega}\partial
_{n}\Phi_{i}\left\vert \nabla f_{i}\right\vert ^{2}=0$ using (\ref{3.3.3}).
The last boundary term is treated as follows. Using $\partial_{n}f_{i}%
=\frac{1}{2}\partial_{n}\Phi_{i}f_{i}$, $\left\vert \Delta\Phi_{i}\right\vert
\leq\frac{cs}{\Gamma}$ and (\ref{3.3.3}), we have
\[%
\begin{array}
[c]{ll}%
\displaystyle\sum_{i=1,..,2d}\displaystyle\int_{\partial\Omega}\partial
_{n}f_{i}\Delta\Phi_{i}f_{i} & =\displaystyle\sum_{i=1,..,2d}\displaystyle\int
_{\partial\Omega}\frac{1}{2}\partial_{n}\Phi_{i}\Delta\Phi_{i}\left\vert
f_{i}\right\vert ^{2}\\
& \leq\displaystyle\frac{cs}{\Gamma}\displaystyle\sum_{i=1,..,d}\int
_{\partial\Omega}\left\vert \partial_{n}\Phi_{i}\right\vert \left\vert
f_{i}\right\vert ^{2}=\displaystyle\frac{cs}{\Gamma}\sum_{i=1,..,d}%
\displaystyle\int_{\partial\Omega}\left(  -\partial_{n}\Phi_{i}\right)
\left\vert f_{i}\right\vert ^{2}%
\end{array}
\]
since $\partial_{n}\psi_{i}\leq0$ and, by an integration by parts
\[%
\begin{array}
[c]{ll}%
\displaystyle\int_{\partial\Omega}\left(  -\partial_{n}\Phi_{i}\right)
\left\vert f_{i}\right\vert ^{2} & =-2\displaystyle\int_{\Omega}\nabla
f_{i}\cdot\nabla\Phi_{i}f_{i}-\displaystyle\int_{\Omega}\Delta\Phi
_{i}\left\vert f_{i}\right\vert ^{2}\\
& \leq\displaystyle\int_{\Omega}\left\vert \nabla f_{i}\right\vert
^{2}+\displaystyle\int_{\Omega}\left\vert \nabla\Phi_{i}\right\vert
^{2}\left\vert f_{i}\right\vert ^{2}+\displaystyle\frac{cs}{h}\left\Vert
\boldsymbol{f}\right\Vert ^{2}%
\end{array}
\]
using Cauchy-Schwarz and $\left\vert \Delta\Phi_{i}\right\vert \leq\frac
{cs}{\Gamma}\leq\frac{cs}{h}$, which implies that
\[
\sum_{i=1,..,2d}\int_{\partial\Omega}\partial_{n}f_{i}\Delta\Phi_{i}f_{i}%
\leq\frac{cs}{\Gamma}\sum_{i=1,..,d}\int_{\Omega}\left\vert \nabla
f_{i}\right\vert ^{2}+\frac{c^{2}s^{2}}{h^{2}}\left\Vert \boldsymbol{f}%
\right\Vert ^{2}+\frac{cs}{\Gamma}\sum_{i=1,..,d}\int_{\Omega}\left\vert
\nabla\Phi_{i}\right\vert ^{2}\left\vert f_{i}\right\vert ^{2}\text{ .}%
\]
One can conclude for the contribution of the boundary terms that for any
$s\in\left(  0,1\right]  $%
\begin{equation}
\text{Boundary terms}\leq\frac{cs}{\Gamma}\sum_{i=1,..,d}\int_{\Omega
}\left\vert \nabla f_{i}\right\vert ^{2}+\frac{c^{2}}{h^{2}}\left\Vert
\boldsymbol{f}\right\Vert ^{2}+\frac{cs}{\Gamma}\sum_{i=1,..,d}\int_{\Omega
}\left\vert \nabla\Phi_{i}\right\vert ^{2}\left\vert f_{i}\right\vert
^{2}\text{ .} \tag{3.5.4}\label{3.5.4}%
\end{equation}

\bigskip

Consequently, from (\ref{3.5.1})-(\ref{3.5.2})-(\ref{3.5.3})-(\ref{3.5.4}), we
obtain that for any $h\in\left(  0,1\right]  $ and any $s\in\left(
0,1\right]  $
\[%
\begin{array}
[c]{ll}%
\left\langle \mathcal{S}^{\prime}\boldsymbol{f},\boldsymbol{f}\right\rangle
+2\left\langle \mathcal{S}\boldsymbol{f},\mathcal{A}\boldsymbol{f}%
\right\rangle  & \leq\displaystyle\frac{cs}{\Gamma}\sum_{i=1,..,2d}%
\int_{\Omega}\left\vert \nabla f_{i}\right\vert ^{2}+\displaystyle\frac{c^{2}%
}{h^{2}}\left\Vert \boldsymbol{f}\right\Vert ^{2}\\
& \quad-\displaystyle\frac{2}{\Gamma}\sum_{i=1,..,2d}\displaystyle\int
_{\Omega}\left(  \eta_{i}+\frac{1}{4}\left\vert \nabla\Phi_{i}\right\vert
^{2}+\frac{s}{4}\nabla\Phi_{i}\nabla^{2}\varphi_{i}\nabla\Phi_{i}\right)
\left\vert f_{i}\right\vert ^{2}\\
& \quad+\displaystyle\frac{cs}{\Gamma}\displaystyle\sum_{i=1,..,d}%
\displaystyle\int_{\Omega}\left\vert \nabla\Phi_{i}\right\vert ^{2}\left\vert
f_{i}\right\vert ^{2}%
\end{array}
\]
which gives that for any $s\in\left(  0,1\right]  $ sufficiently small,
\begin{equation}%
\begin{array}
[c]{ll}%
\left\langle \mathcal{S}^{\prime}\boldsymbol{f},\boldsymbol{f}\right\rangle
+2\left\langle \mathcal{S}\boldsymbol{f},\mathcal{A}\boldsymbol{f}%
\right\rangle  & \leq\displaystyle\frac{Cs}{\Gamma}\sum_{i=1,..,2d}%
\displaystyle\int_{\Omega}\left\vert \nabla f_{i}\right\vert ^{2}+\frac
{C}{h^{2}}\left\Vert \boldsymbol{f}\right\Vert ^{2}\\
& -\displaystyle\frac{2}{\Gamma}\sum_{i=1,..,2d}\displaystyle\int_{\Omega
}\left(  \eta_{i}+\frac{1}{8}\left\vert \nabla\Phi_{i}\right\vert ^{2}\right)
\left\vert f_{i}\right\vert ^{2}\text{ .}%
\end{array}
\tag{3.5.5}\label{3.5.5}%
\end{equation}
Indeed, $-\frac{s}{4}\nabla\Phi_{i}\nabla^{2}\varphi_{i}\nabla\Phi_{i}%
\leq\frac{s}{4}\left\vert \nabla^{2}\varphi_{i}\right\vert \left\vert
\nabla\Phi_{i}\right\vert ^{2}\leq cs\left\vert \nabla\Phi_{i}\right\vert
^{2}$.

\bigskip

It remains to prove that
\begin{equation}
-\frac{2}{\Gamma}\sum_{i=1,..,2d}\int_{\Omega}\left(  \eta_{i}+\frac{1}%
{8}\left\vert \nabla\Phi_{i}\right\vert ^{2}\right)  \left\vert f_{i}%
\right\vert ^{2}\leq C\left\Vert \boldsymbol{f}\right\Vert ^{2}+\frac
{2-s/c}{\Gamma}\sum_{i=1,..,2d}\int_{\Omega}\left(  -\eta_{i}\right)
\left\vert f_{i}\right\vert ^{2}\text{ .} \tag{3.5.6}\label{3.5.6}%
\end{equation}

\bigskip

By Proposition 5 $\left(  i\right)  $ and $\left(  ii\right)  $, $\left\vert
\varphi_{i,1}\right\vert \leq\frac{c}{2}\left\vert \nabla\varphi
_{i,1}\right\vert ^{2}$ in $\mathcal{B}_{i}\cup\mathcal{D}_{i}$. This implies
that for any $i\in\left\{  1,\cdot\cdot,d\right\}  $
\[
-\left\vert \nabla\Phi_{i}\right\vert ^{2}=-\frac{s^{2}}{\Gamma^{2}}\left\vert
\nabla\varphi_{i,1}\right\vert ^{2}\leq-\frac{2s^{2}}{c\Gamma^{2}}\left\vert
\varphi_{i,1}\right\vert =\frac{4s}{c}\left(  -\frac{s}{2\Gamma^{2}}\left\vert
\varphi_{i,1}\right\vert \right)  \leq\frac{4s}{c}\eta_{i}\text{ .}%
\]
Therefore, we get that for any $i\in\left\{  1,\cdot\cdot,d\right\}  $
\[
-\frac{1}{4}\int_{\mathcal{B}_{i}\cup\mathcal{D}_{i}}\left\vert \nabla\Phi
_{i}\right\vert ^{2}\left\vert f_{i}\right\vert ^{2}\leq\frac{s}{c}%
\int_{\mathcal{B}_{i}\cup\mathcal{D}_{i}}\eta_{i}\left\vert f_{i}\right\vert
^{2}%
\]
which yields%
\[
-\displaystyle\frac{2}{\Gamma}\sum_{i=1,..,d}\displaystyle\int_{\mathcal{B}%
_{i}\cup\mathcal{D}_{i}}\left(  \eta_{i}+\frac{1}{8}\left\vert \nabla\Phi
_{i}\right\vert ^{2}\right)  \left\vert f_{i}\right\vert ^{2}\leq
\displaystyle\frac{2-s/c}{\Gamma}\displaystyle\sum_{i=1,..,d}\int
_{\mathcal{B}_{i}\cup\mathcal{D}_{i}}\left(  -\eta_{i}\right)  \left\vert
f_{i}\right\vert ^{2}\text{ .}%
\]
By Proposition 5 $\left(  iv\right)  $, there is $c_{3}>0$ such that
$\varphi_{i,1}-\varphi_{j,1}\leq-c_{3}$ in $\mathcal{C}_{i}$ for some $j\neq
i$. Therefore, using $\left\vert \eta_{i}+\frac{1}{8}\left\vert \nabla\Phi
_{i}\right\vert ^{2}\right\vert \leq\frac{c}{\Gamma^{2}}$ and $\left\vert
f_{i}\right\vert ^{2}=e^{s\left(  \varphi_{i,1}-\varphi_{j,1}\right)  \frac
{1}{\Gamma}}\left\vert f_{j}\right\vert ^{2}$, it holds
\[
-\displaystyle\frac{2}{\Gamma}\sum_{i=1,..,d}\displaystyle\int_{\mathcal{C}%
_{i}}\left(  \eta_{i}+\frac{1}{8}\left\vert \nabla\Phi_{i}\right\vert
^{2}\right)  \left\vert f_{i}\right\vert ^{2}\leq\displaystyle\frac{2c}%
{\Gamma^{3}}e^{-c_{3}\frac{s}{\Gamma}}\left(  \sum_{i=1,..,d}\displaystyle\int
_{%
{\textstyle\bigcup\limits_{j\neq i}}
\Theta_{p_{j}}}\left\vert f_{j}\right\vert ^{2}\right)  \leq C_{s}\left\Vert
\boldsymbol{f}\right\Vert ^{2}\text{ .}%
\]
By Proposition 5 $\left(  v\right)  $, $\left\vert \varphi_{i,2}\right\vert
\leq c_{5}\left\vert \nabla\varphi_{i,2}\right\vert ^{2}$ in a neighborhood
$\vartheta$ of $\partial\Omega$ and similarly one can deduce that,
\[
-\displaystyle\frac{2}{\Gamma}\sum_{i=d+1,..,2d}\displaystyle\int_{\vartheta
}\left(  \eta_{i}+\frac{1}{8}\left\vert \nabla\Phi_{i}\right\vert ^{2}\right)
\left\vert f_{i}\right\vert ^{2}\leq\displaystyle\frac{2-s/c}{\Gamma}%
\sum_{i=d+1,..,2d}\displaystyle\int_{\vartheta}\left(  -\eta_{i}\right)
\left\vert f_{i}\right\vert ^{2}\text{ .}%
\]
By Proposition 5 $\left(  vi\right)  $, there is $c_{6}>0$ such that
$\varphi_{i,2}-\varphi_{i,1}\leq-c_{6}$ outside the neighborhood $\vartheta$
of $\partial\Omega$ which implies%
\[
-\displaystyle\frac{2}{\Gamma}\sum_{i=d+1,..,2d}\displaystyle\int
_{\Omega\left\backslash \vartheta\right.  }\left(  \eta_{i}+\frac{1}%
{8}\left\vert \nabla\Phi_{i}\right\vert ^{2}\right)  \left\vert f_{i}%
\right\vert ^{2}\leq\displaystyle\frac{2c}{\Gamma^{3}}e^{-c_{6}\frac{s}%
{\Gamma}}\sum_{i=1,..,d}\displaystyle\int_{\Omega\left\backslash
\vartheta\right.  }\left\vert f_{i}\right\vert ^{2}\leq C_{s}\left\Vert
\boldsymbol{f}\right\Vert ^{2}\text{ .}%
\]
This completes the proof of (\ref{3.5.6}).

\bigskip

Consequently, by (\ref{3.5.5}) and (\ref{3.5.6}) one can conclude that for any
$h\in\left(  0,1\right]  $ and any $s\in\left(  0,1\right]  $ sufficiently
small,
\[
\left\langle \mathcal{S}^{\prime}\boldsymbol{f},\boldsymbol{f}\right\rangle
+2\left\langle \mathcal{S}\boldsymbol{f},\mathcal{A}\boldsymbol{f}%
\right\rangle \leq\frac{Cs}{\Gamma}\sum_{i=1,..,2d}\int_{\Omega}\left\vert
\nabla f_{i}\right\vert ^{2}+\frac{C_{s}}{h^{2}}\left\Vert \boldsymbol{f}%
\right\Vert ^{2}+\frac{2-s/c}{\Gamma}\sum_{i=1,..,2d}\int_{\Omega}\left(
-\eta_{i}\right)  \left\vert f_{i}\right\vert ^{2}%
\]
which implies%
\[
\left\langle \mathcal{S}^{\prime}\boldsymbol{f},\boldsymbol{f}\right\rangle
+2\left\langle \mathcal{S}\boldsymbol{f},\mathcal{A}\boldsymbol{f}%
\right\rangle \leq\frac{1+C_{0}}{\Gamma}\left\langle \mathcal{S}%
\boldsymbol{f},\boldsymbol{f}\right\rangle +\frac{C}{h^{2}}\left\Vert
\boldsymbol{f}\right\Vert ^{2}%
\]
with $C_{0}\in\left(  0,1\right)  $ and $C>0$. Finally, the system
(\ref{3.4.1}) of ordinary differential inequalities becomes%
\[
\left\{
\begin{array}
[c]{ll}%
\left\vert \displaystyle\frac{1}{2}\frac{d}{dt}\left\Vert \boldsymbol{f}%
\right\Vert ^{2}+\mathbf{N}\left(  t\right)  \left\Vert \boldsymbol{f}%
\right\Vert ^{2}\right\vert \leq\left\Vert a\right\Vert _{\infty}\left\Vert
\boldsymbol{f}\right\Vert ^{2}\text{ ,} & \\
\mathbf{N}^{\prime}\left(  t\right)  \leq\displaystyle\frac{1+C_{0}}{\Gamma
}\mathbf{N}\left(  t\right)  +\left\Vert a\right\Vert _{\infty}^{2}%
+\displaystyle\frac{C}{h^{2}}\text{ .} &
\end{array}
\right.
\]

\bigskip

\subsection{Step 6: Solve ODE}

\bigskip

Let $h\in\left(  0,1\right]  $\ and $\ell>1$\ such that $\ell h<$min$\left(
1/2,T/4\right)  $. Applying Proposition 3 with $t_{3}=T$, $t_{2}=T-\ell h$,
and $t_{1}=T-2\ell h$, we obtain that
\[
y\left(  T-\ell h\right)  ^{1+M_{\ell}}\leq y\left(  T\right)  y\left(
T-2\ell h\right)  ^{M_{\ell}}e^{D_{\ell}}%
\]
where $D_{\ell}=2M_{\ell}\left(  F_{2}\left(  2\ell h\right)  ^{2}%
+F_{1}\left(  2\ell h\right)  \right)  $, $M_{\ell}=\frac{\left(
\ell+1\right)  ^{C_{0}}-1}{1-\left(  \frac{\ell+1}{2\ell+1}\right)  ^{C_{0}}%
}\leq\frac{\left(  \ell+1\right)  ^{C_{0}}}{1-\left(  \frac{2}{3}\right)
^{C_{0}}}$\ if $C_{0}>0$.

\bigskip

From now, $y\left(  t\right)  =\left\Vert \boldsymbol{f}\left(  \cdot
,t\right)  \right\Vert ^{2}$, $N$ is the frequency function $\mathbf{N}$,
$F_{1}=\left\Vert a\right\Vert _{\infty}$ and $F_{2}=\left\Vert a\right\Vert
_{\infty}^{2}+\frac{C}{h^{2}}$: We have by the above Proposition 3 and Step 5,%
\begin{equation}
\left(  \left\Vert \boldsymbol{f}\left(  \cdot,T-\ell h\right)  \right\Vert
^{2}\right)  ^{1+M_{\ell}}\leq\left\Vert \boldsymbol{f}\left(  \cdot,T\right)
\right\Vert ^{2}\left(  \left\Vert \boldsymbol{f}\left(  \cdot,T-2\ell
h\right)  \right\Vert ^{2}\right)  ^{M_{\ell}}K_{\ell}\text{ } \tag{3.6.1}%
\label{3.6.1}%
\end{equation}
where $K_{\ell}=e^{D_{\ell}}$ with $D_{\ell}=2M_{\ell}\left(  \left(
\left\Vert a\right\Vert _{\infty}^{2}+\frac{C}{h^{2}}\right)  \left(  2\ell
h\right)  ^{2}+\left\Vert a\right\Vert _{\infty}\left(  2\ell h\right)
\right)  $. Notice that when $\left\Vert a\right\Vert _{\infty}^{2/3}h<1$,
then the following upper bound for $K_{\ell}$ holds
\begin{equation}
K_{\ell}\leq e^{C_{\ell}\left(  1+\left\Vert a\right\Vert _{\infty}%
^{2/3}\right)  }\text{ .} \tag{3.6.2}\label{3.6.2}%
\end{equation}
Indeed, $D_{\ell}\leq2M_{\ell}\left(  1+4C\ell^{2}+2\left\Vert a\right\Vert
_{\infty}^{2}\left(  2\ell h\right)  ^{2}\right)  $ and $h^{2}\left\Vert
a\right\Vert _{\infty}^{2}=\left\Vert a\right\Vert _{\infty}^{2/3}\left(
\left\Vert a\right\Vert _{\infty}^{2/3}h\right)  ^{2}\leq\left\Vert
a\right\Vert _{\infty}^{2/3}$.

\bigskip

\subsection{Step 7: Make appear $\omega$}

\bigskip

It is well-known that for any $0\leq t_{1}\leq t_{2}\leq T$,
\begin{equation}
\left\Vert u\left(  \cdot,t_{2}\right)  \right\Vert _{L^{2}\left(
\Omega\right)  }\leq e^{\left(  t_{2}-t_{1}\right)  \left\Vert a\right\Vert
_{\infty}}\left\Vert u\left(  \cdot,t_{1}\right)  \right\Vert _{L^{2}\left(
\Omega\right)  } \tag{3.7.1}\label{3.7.1}%
\end{equation}
where $\left\Vert a\right\Vert _{\infty}=\left\Vert a\right\Vert _{L^{\infty
}\left(  \Omega\times\left(  0,T\right)  \right)  }$.

\bigskip

Observe that
\[
\left\Vert f_{1}\right\Vert _{L^{2}\left(  \Omega\right)  }^{2}\leq\left\Vert
\boldsymbol{f}\right\Vert ^{2}\leq2\sum_{i=1,..,d}\left\Vert f_{i}\right\Vert
_{L^{2}\left(  \Omega\right)  }^{2}%
\]
since $\varphi_{i,2}\leq\varphi_{i,1}$ on $\Omega$. Therefore, (\ref{3.6.1})
becomes
\begin{equation}%
\begin{array}
[c]{ll}%
\left(  \left\Vert f_{1}\left(  \cdot,T-\ell h\right)  \right\Vert
_{L^{2}\left(  \Omega\right)  }^{2}\right)  ^{1+M_{\ell}} & \leq
2\displaystyle\sum_{i=1,..,d}\left\Vert f_{i}\left(  \cdot,T\right)
\right\Vert _{L^{2}\left(  \Omega\right)  }^{2}\\
& \quad\times\left(  2\displaystyle\sum_{i=1,..,d}\left\Vert f_{i}\left(
\cdot,T-2\ell h\right)  \right\Vert _{L^{2}\left(  \Omega\right)  }%
^{2}\right)  ^{M_{\ell}}K_{\ell}\text{ .}%
\end{array}
\tag{3.7.2}\label{3.7.2}%
\end{equation}
First, notice that from (\ref{3.7.1}), using $\Phi_{i}\leq0$,
\begin{equation}
\left\Vert f_{i}\left(  \cdot,T-2\ell h\right)  \right\Vert _{L^{2}\left(
\Omega\right)  }^{2}\leq e^{2T\left\Vert a\right\Vert _{\infty}}\int_{\Omega
}\left\vert u\left(  \cdot,0\right)  \right\vert ^{2}\text{ .} \tag{3.7.3}%
\label{3.7.3}%
\end{equation}
Second, we make appear $\omega_{i,r}=\left\{  x;\left\vert x-p_{i}\right\vert
<r\right\}  \subset\omega$ from $\left\Vert f_{i}\left(  \cdot,T\right)
\right\Vert _{L^{2}\left(  \Omega\right)  }^{2}$ as follows:%
\begin{equation}%
\begin{array}
[c]{ll}%
\left\Vert f_{i}\left(  \cdot,T\right)  \right\Vert _{L^{2}\left(
\Omega\right)  }^{2} & =\displaystyle\int_{\omega_{i,r}}\left\vert u\left(
\cdot,T\right)  \right\vert ^{2}e^{\frac{s}{h}\varphi_{i,1}}+\displaystyle\int
_{\Omega\left\backslash \omega_{i,r}\right.  }\left\vert u\left(
\cdot,T\right)  \right\vert ^{2}e^{\frac{s}{h}\varphi_{i,1}}\\
& \leq\displaystyle\int_{\omega}\left\vert u\left(  \cdot,T\right)
\right\vert ^{2}+e^{-\frac{s\mu}{h}}\displaystyle e^{2T\left\Vert a\right\Vert
_{\infty}}\int_{\Omega}\left\vert u\left(  \cdot,0\right)  \right\vert ^{2}%
\end{array}
\tag{3.7.4}\label{3.7.4}%
\end{equation}
because on $\Omega\left\backslash \omega_{i,r}\right.  $, $\varphi_{i,1}%
\leq-\mu$ for some $\mu>0$ and we used (\ref{3.7.1}). Third, from
(\ref{3.7.1}) with $\ell h<$min$\left(  1/2,T/4\right)  $ and $-\varphi
_{1,1}\leq c$ it holds%
\begin{equation}%
\begin{array}
[c]{ll}%
\displaystyle\int_{\Omega}\left\vert u\left(  \cdot,T\right)  \right\vert ^{2}
& \leq\displaystyle e^{2\ell h\left\Vert a\right\Vert _{\infty}}\int_{\Omega
}\left\vert u\left(  \cdot,T-\ell h\right)  \right\vert ^{2}e^{\frac
{s}{\left(  \ell+1\right)  h}\varphi_{1,1}}e^{-\frac{s}{\left(  \ell+1\right)
h}\varphi_{1,1}}\\
& \leq e^{T\left\Vert a\right\Vert _{\infty}}e^{\frac{sc}{\left(
\ell+1\right)  h}}\left\Vert f_{1}\left(  \cdot,T-\ell h\right)  \right\Vert
_{L^{2}\left(  \Omega\right)  }^{2}\text{ .}%
\end{array}
\tag{3.7.5}\label{3.7.5}%
\end{equation}
Combining the above four facts (\ref{3.7.2}), (\ref{3.7.3}), (\ref{3.7.4}) and
(\ref{3.7.5}), we can deduce that%
\[%
\begin{array}
[c]{ll}
& \quad\left(  \displaystyle\int_{\Omega}\left\vert u\left(  \cdot,T\right)
\right\vert ^{2}\right)  ^{1+M_{\ell}}\\
& \leq e^{\frac{sc\left(  1+M_{\ell}\right)  }{\left(  \ell+1\right)  h}%
}e^{T\left\Vert a\right\Vert _{\infty}\left(  1+M_{\ell}\right)  }\left(
\left\Vert f_{1}\left(  \cdot,T-\ell h\right)  \right\Vert _{L^{2}\left(
\Omega\right)  }^{2}\right)  ^{1+M_{\ell}}\\
& \leq e^{\frac{sc\left(  1+M_{\ell}\right)  }{\left(  \ell+1\right)  h}%
}e^{T\left\Vert a\right\Vert _{\infty}\left(  1+M_{\ell}\right)  }\left(
2\displaystyle\sum_{i=1,..,d}\left\Vert f_{i}\left(  \cdot,T-2\ell h\right)
\right\Vert _{L^{2}\left(  \Omega\right)  }^{2}\right)  ^{M_{\ell}}K_{\ell}\\
& \quad\times\left(  2\displaystyle\sum_{i=1,..,d}\left\Vert f_{i}\left(
\cdot,T\right)  \right\Vert _{L^{2}\left(  \Omega\right)  }^{2}\right) \\
& \leq e^{\frac{sc\left(  1+M_{\ell}\right)  }{\left(  \ell+1\right)  h}%
}e^{T\left\Vert a\right\Vert _{\infty}\left(  1+M_{\ell}\right)  }\left(
2de^{2T\left\Vert a\right\Vert _{\infty}}\displaystyle\int_{\Omega}\left\vert
u\left(  \cdot,0\right)  \right\vert ^{2}\right)  ^{M_{\ell}}K_{\ell}\\
& \quad\times2d\left(  \displaystyle\int_{\omega}\left\vert u\left(
\cdot,T\right)  \right\vert ^{2}+e^{-\frac{s\mu}{h}}%
\displaystyle e^{2T\left\Vert a\right\Vert _{\infty}}\int_{\Omega}\left\vert
u\left(  \cdot,0\right)  \right\vert ^{2}\right)  \text{ .}%
\end{array}
\]

\bigskip

We will choose $\ell>1$ large enough in order that $\frac{sc\left(  1+M_{\ell
}\right)  }{\left(  \ell+1\right)  h}-\frac{s\mu}{h}\leq-\frac{s\mu}{2h}$ that
is $\frac{c\left(  1+M_{\ell}\right)  }{\left(  \ell+1\right)  }\leq\frac{\mu
}{2}$. This is possible because $M_{\ell}\leq\frac{\left(  \ell+1\right)
^{C_{0}}}{1-\left(  \frac{2}{3}\right)  ^{C_{0}}}$ with $C_{0}\in\left(
0,1\right)  $. Therefore, combining with the upper bound for $K_{\ell}$ (see
(\ref{3.6.2})), there are $M>0$ and $c>0$, such that for any $h>0$\ satisfying
$\ell h<$min$\left(  1/2,T/4\right)  $ and $\left\Vert a\right\Vert _{\infty
}^{2/3}h<1$, we have%
\[%
\begin{array}
[c]{ll}%
\left(  \displaystyle\int_{\Omega}\left\vert u\left(  \cdot,T\right)
\right\vert ^{2}\right)  ^{1+M} & \leq e^{c\left(  1+T\left\Vert a\right\Vert
_{\infty}+\left\Vert a\right\Vert _{\infty}^{2/3}\right)  }\left(
\displaystyle\int_{\Omega}\left\vert u\left(  \cdot,0\right)  \right\vert
^{2}\right)  ^{M}\\
& \quad\times\left(  e^{\frac{s\mu}{2h}}\displaystyle\int_{\omega}\left\vert
u\left(  \cdot,T\right)  \right\vert ^{2}+e^{-\frac{s\mu}{2h}}%
\displaystyle\int_{\Omega}\left\vert u\left(  \cdot,0\right)  \right\vert
^{2}\right)  \text{ .}%
\end{array}
\]
On the other hand, using (\ref{3.7.1}), for any $h\geq$min$\left(  1/\left(
2\ell\right)  ,T/\left(  4\ell\right)  \right)  $,
\[
\int_{\Omega}\left\vert u\left(  \cdot,T\right)  \right\vert ^{2}\leq
e^{2T\left\Vert a\right\Vert _{\infty}}\int_{\Omega}\left\vert u\left(
\cdot,0\right)  \right\vert ^{2}e^{-\frac{s\mu}{2h}}e^{\frac{s\mu}{2}\left(
2\ell+\frac{4\ell}{T}\right)  }\text{ ,}%
\]
and for any $h$ such that $1\leq\left\Vert a\right\Vert _{\infty}^{2/3}h$,
\[
\int_{\Omega}\left\vert u\left(  \cdot,T\right)  \right\vert ^{2}\leq
e^{2T\left\Vert a\right\Vert _{\infty}}\int_{\Omega}\left\vert u\left(
\cdot,0\right)  \right\vert ^{2}e^{-\frac{s\mu}{2h}}e^{\frac{s\mu}%
{2}\left\Vert a\right\Vert _{\infty}^{2/3}}\text{ .}%
\]
Consequently, one can conclude that for any $h>0$, it holds%
\[%
\begin{array}
[c]{ll}%
\left(  \displaystyle\int_{\Omega}\left\vert u\left(  \cdot,T\right)
\right\vert ^{2}\right)  ^{1+M} & \leq e^{c\left(  1+\frac{1}{T}+T\left\Vert
a\right\Vert _{\infty}+\left\Vert a\right\Vert _{\infty}^{2/3}\right)
}\left(  \displaystyle\int_{\Omega}\left\vert u\left(  \cdot,0\right)
\right\vert ^{2}\right)  ^{M}\\
& \quad\times\left(  e^{\frac{c}{h}}\displaystyle\int_{\omega}\left\vert
u\left(  \cdot,T\right)  \right\vert ^{2}+e^{-\frac{1}{h}}\displaystyle\int
_{\Omega}\left\vert u\left(  \cdot,0\right)  \right\vert ^{2}\right)  \text{
.}%
\end{array}
\]
Now, choose $h>0$ such that
\[
e^{-\frac{1}{h}}e^{c\left(  1+\frac{1}{T}+T\left\Vert a\right\Vert _{\infty
}+\left\Vert a\right\Vert _{\infty}^{2/3}\right)  }\left(  e^{2T\left\Vert
a\right\Vert _{\infty}}\int_{\Omega}\left\vert u\left(  \cdot,0\right)
\right\vert ^{2}\right)  ^{1+M}=\frac{1}{2}\left(  \int_{\Omega}\left\vert
u\left(  \cdot,T\right)  \right\vert ^{2}\right)  ^{1+M}\text{ ,}%
\]
we obtain the desired estimate for some $M_{1}>1$ and $c_{1}>0$%
\[
\left(  \int_{\Omega}\left\vert u\left(  \cdot,T\right)  \right\vert
^{2}\right)  ^{1+M_{1}}\leq e^{c_{1}\left(  1+\frac{1}{T}+T\left\Vert
a\right\Vert _{\infty}+\left\Vert a\right\Vert _{\infty}^{2/3}\right)  }%
\int_{\omega}\left\vert u\left(  \cdot,T\right)  \right\vert ^{2}\left(
\int_{\Omega}\left\vert u\left(  \cdot,0\right)  \right\vert ^{2}\right)
^{M_{1}}\text{ .}%
\]
This completes the proof.

\bigskip

\bigskip

\bigskip

\bigskip

\bigskip


\bigskip

\bigskip

\bigskip

\begin{thebibliography}{9999}                                                                                             %


\bibitem[AEWZ]{AEWZ}J. Apraiz, L. Escauriaza, G. Wang and C. Zhang,
Observability inequalities and measurable sets, J. Eur. Math. Soc. (JEMS) 16
(11) (2014), 2433--2475.

\bibitem[BP]{BP}C. Bardos and K.D. Phung, Observation estimate for kinetic
transport equations by diffusion approximation, C. R. Math. Acad. Sci. Paris
355 (2017), no.6, 640--664.

\bibitem[BuP]{BuP}R. Buffe and K.D. Phung, A spectral inequality for
degenerate operators and applications, C. R. Math. Acad. Sci. Paris 356
(11-12) (2018) 1131--1155.

\bibitem[BM]{BM}N. Burq and I. Moyano, Propagation of smallness and control
for heat equations, ArXiv:1912.07402.

\bibitem[C]{C}J.-M. Coron, Control and Nonlinearity, vol. 136 of Mathematical
Surveys and Monographs, AMS, Providence, RI, 2007.

\bibitem[DZZ]{DZZ}T. Duyckaerts, X. Zhang and E. Zuazua, On the optimality of
the observability inequalities for parabolic and hyperbolic systems with
potentials, Ann. Inst. H. Poincar\'{e} Anal. Non Lin\'{e}aire 25 (2008) 1--41.

\bibitem[EFV]{EFV}L. Escauriaza, F.J. Fern\'{a}ndez and S. Vessella, Doubling
properties of caloric functions, Appl. Anal. 85 (2006) 205--223.

\bibitem[FGGP]{FGGP}E. Fern\'{a}ndez-Cara, M. Gonz\'{a}lez-Burgos, S. Guerrero
and J.-P. Puel, Null controllability of the heat equation with boundary
Fourier conditions: the linear case, ESAIM Control Optim. Calc. Var. 12(3)
(2006) 442--465.

\bibitem[FV]{FV}E. Francini and S. Vessella, Carleman estimates for the
parabolic transmission problem and H\"{o}lder propagation of smallness across
an interface, J. Differential Equations 265 (2018) 2375--2430.

\bibitem[FI]{FI}A.V. Fursikov and O.Y. Imanuvilov, Controllability of
evolution equations. Lecture Notes Series, 34. Seoul National University,
Research Institute of Mathematics, Global Analysis Research Center, Seoul, 1996.

\bibitem[K]{K}K. Kurata, On a backward estimate for solutions of parabolic
differential equations and its application to unique continuation. Spectral
and scattering theory and applications, 247--257, Adv. Stud. Pure Math., 23,
Math. Soc. Japan, Tokyo, 1994.

\bibitem[LRL]{LRL}J. Le Rousseau and G. Lebeau, On Carleman estimates for
elliptic and parabolic operators. Applications to unique continuation and
control of parabolic equations, ESAIM Control Optim. Calc. Var., 18 (2012), 712--747.

\bibitem[LR]{LR}G. Lebeau and L. Robbiano, Contr\^{o}le exact de
l'\'{e}quation de la chaleur, Communications in Partial Differential Equation,
20 (1995) 335--356.

\bibitem[M]{M}Y. Matsumoto, An introduction to Morse Theory, American
Mathematical Society, 2002

\bibitem[P]{P}K.D. Phung, Carleman commutator approach in logarithmic
convexity for parabolic equations, Mathematical Control and Related Fields 8
(3-4) (2018) 899--933.

\bibitem[PW]{PW}K. D. Phung and G. Wang, An observability estimate for
parabolic equations from a measurable set in time and its applications, J.
Eur. Math. Soc. (JEMS) 15 (2) (2013) 681--703.

\bibitem[PWZ]{PWZ}K. D. Phung, L. Wang and C. Zhang, Bang-bang property for
time optimal control of semilinear heat equation, Ann. Inst. H. Poincar\'{e}
Anal. Non Lin\'{e}aire 31 (3) (2014) 477--499.

\bibitem[TW]{TW}M. Tucsnak and G. Weiss, Observation and control for operator
semigroups, Birkh\"{a}user Advanced Texts: Basler Lehrb\"{u}cher,
Birkh\"{a}user Verlag, Basel, 2009.

\bibitem[WW]{WW}G. Wang and L. Wang, The Carleman inequality and its
application to periodic optimal control governed by semilinear parabolic
differential equations. J. Optim. Theory Appl. 118 (2003), no. 2, 429--461.
\end{thebibliography}
\end{document}